\documentclass[12pt,reqno]{article}
\usepackage{amssymb,amscd,amsmath,amsthm,color}
\newtheorem{theorem}{Theorem}[section]
\newtheorem{lemma}[theorem]{Lemma}
\newtheorem{cor}[theorem]{Corollary}

\usepackage{enumerate,amssymb}
\theoremstyle{definition}

\theoremstyle{remark}

\numberwithin{equation}{section}

\def\bC{\mathbb{C}}

\def\bM{\mathbb{M}}

\begin{document}
\baselineskip=15pt

\title{Numerical range and positive block matrices}

\author{ Jean-Christophe Bourin{\footnote{Funded by the ANR Projet (No.\ ANR-19-CE40-0002) and by the French Investissements
 d'Avenir program, project ISITE-BFC (contract ANR-15-IDEX-03).
}}\,  and Eun-Young Lee{\footnote{This research was supported by
Basic Science Research Program through the National Research
Foundation of Korea (NRF) funded by the Ministry of
Education (NRF-2018R1D1A3B07043682)}  }  }

\date{ }

\maketitle

\vskip 10pt\noindent
{\small 
{\bf Abstract.} We obtain several norm and eigenvalue inequalities for positive matrices partitioned into four blocks. The results involve the numerical range $W(X$)  of the off-diagonal block $X$, especially the distance $d$ from $0$ to $W(X)$.
 A special consequence is an estimate,
$$ 
{\mathrm{diam}}\, W\left(\begin{bmatrix} A &X \\
X^* &B\end{bmatrix}\right)
- {\mathrm{diam}}\, W\left(\frac{A+B}{2}\right) \ge 2d,
$$
between the diameters of the numerical ranges for the full matrix and its partial trace.


}

\section{The width of the numerical range}

Let $\bM_n$ denote the space of $n$-by-$n$ matrices, and let $\langle h, h' \rangle$  be the canonical inner product of $\bC^n$, linear in the first variable.  The numerical range of $X\in\bM_n$ is defined as 
$$
W(X) =\{ \langle h, Xh\rangle \ : \ \| h\| =1\}.
$$
The Hausdorff-Toeplitz theorem states that $W(X)$ is a compact convex set containing the spectrum of $X$. In case of a normal matrix, the numerical range is precisely the convex hull of the spectrum. The symbol $\|\cdot\|$ will also denote any symmetric norm on $\bM_{2n}$. Such a norm is also called a unitarily invariant norm. It  sastifies the unitary invariance property $\| UTV\|=\|T\|$ for all $T\in\bM_{2n}$ and all unitary matrices $U,V\in\bM_{2n}$, and  it induces a symmetric norm on $\bM_n$ in an obvious way, by considering $\bM_n$ as the upper left corner of $\bM_{2n}$ completed with some zero entries.

It has been recently pointed out \cite{BM} that the numerical range plays a  role to estimate  a partitioned positive matrix with
 its partial trace, i.e.,  the sum of the diagonal blocks.  In Matrix Analysis, positive matrices partitioned into four blocks  are a fundamental tool and these matrices are also of basic importance in applications, especially in Quantum Information Theory. The main theorem of \cite{BM} reads as follows.

\vskip 5pt
\begin{theorem}\label{th1} Let $\begin{bmatrix} A &X \\
X^* &B\end{bmatrix} $ be a positive matrix  partitioned into four blocks
in $\bM_n$. Suppose that $W(X)$ has the width $\omega$. Then, for all symmetric norms,
$$
\left\| \begin{bmatrix} A &X \\
X^* &B\end{bmatrix}\right\| \le \| A+B +\omega I \|.
$$
\end{theorem}

\vskip 5pt
Here $I$ stands for the identity matrix and the width of $W(X)$ is the smallest distance between two parallel straight lines such that the strip between these two lines contains $W(X)$. Hence the partial trace $A+B$ may be used to give an upper bound for the norms of the full block-matrix. This note will provide a lower bound, stated in Section 2, and several consequences.

 Theorem \ref{th1} is   the first  inequality involving the width of the numerical range; classical results  rather deal with the numerical radius, $w(X)=\max\{|z| \, :\, z\in W(X)\}.
$
Our new lower bound will also have an unusual feature as it involves the distance from $0$ to the numerical range, ${\mathrm{dist}}(0,W(X))=\min\{|z| \, :\, z\in W(X)\}.$ For a background on the numerical range we refer to \cite{HJ}, where the term of Field of values is used. Some very interesting inequalities for the numerical radius can be found in \cite{Hol}, \cite{K1}, and in the recent article \cite{Cain}.

In case of  Hermitian off-diagonal blocks, Theorem \ref{th1} holds with $w=0$. More generally, if $X=X=a I +bH$ for some  $a,b\in\bC$ and some Hermitian matrix $H$, we have $\omega=0$ as $W(X)$ is a line segment. This special case of the theorem was first shown by Mhanna \cite{Mha}. In particular, if the  off-diagonal blocks are normal two-by-two matrices, then we can take $\omega=0$. This does not hold any longer for three-by-three normal matrices, a detailed study of this phenomenon is given in \cite{Hay} and \cite{GLRT}.

For Hermitian off-diagonal blocks, a stronger statement than Theorem \ref{th1} with  $w=0$  holds. The following decomposition was shown in \cite[Theorem 2.2]{BLL2}.

\vskip 5pt
\begin{theorem}\label{th2} Let $\begin{bmatrix} A &X \\
X &B\end{bmatrix} $ be a positive matrix  partitioned into four Hermitian blocks 
in $\bM_n$.  Then, for some pair of unitary matrices $U,V\in\bM_{2n}$,
$$
 \begin{bmatrix} A &X \\
X &B\end{bmatrix} = \frac{1}{2}\left\{ U  \begin{bmatrix} A+B &0 \\
0 &0\end{bmatrix} U^* 
+ V \begin{bmatrix} 0 &0\\
0 &A+B\end{bmatrix} V^*
\right\}.
$$
\end{theorem}

\vskip 5pt
 For  decompositions of positive matrices partitioned into a larger number of blocks, see \cite{BL}. We close this  section by recalling some facts on symmetric norms, classical text books such as \cite{Bhatia}, \cite{HJ} and \cite{S} are good references.

A symmetric norm  on $\bM_n$, can be defined
by its restriction to the positive cone $\bM_n^+$. Symmetric norms  on $\bM_n^+$
are characterized by  three properties:
\begin{itemize}
\item[(i)] $\| \lambda A\| = \lambda\|A\|$  for all $A\in \bM_n^+ $ and all  $\lambda\ge 0$,
\item[(ii)] $\| UAU^*\|$ for all $A\in \bM_n^+$ and all unitaries $U\in\bM_n$, 
\item[(iii)] $\|A\| \le \|  A+B\| \le   \|  A\|  + \|  B\|$ for all $A,\,B\in \bM_n^+ $.
\end{itemize}
 Let $\lambda_1^{\downarrow}(A)\ge\cdots\ge \lambda_n^{\downarrow}(A)$ stand for the eigenvalues of $A\in\bM_n^+$ arranged in non-increasing order. Then, the Ky Fan $k$-norms,
$$
\| A\|_{(k)} =\sum_{j=1}^k \lambda_j^{\downarrow}(A)
$$
are symmetric norms, $k=1,\ldots,n$. Thus $\| A\|_{(1)}$ is the operator norm, usually denoted by $\| A\|_{\infty}$ while $\| A\|_{(n)}$ is the trace norm, usually written $\| A\|_1$. For $A,B\in\bM_n^+$, the following conditions are equivalent:
\begin{itemize}
\item[(a)] $\| A\|_{(k)} \le \| B\|_{(k)} $ for all $k=1,\ldots,n$,
\item[(b)] $\| A\| \le \| B\|$ for all symmetric norms,
\item[(c)] The vector of the eigenvalues of $A$ is dominated by a convex combination of permutations of the vector of the eigenvalues of $B$, equivalently, 
$$
A \le \sum_{i=1}^{n+1} \alpha_i U_i BU_i^*
$$
for some  unitary matrices $U_i$ and some  scalars $\alpha_i\ge 0$ such that $\sum_{i=1}^{n+1}\alpha_i=1$.
\end{itemize}
When these conditions hold (especially when explicitly stated as (a)) one says that $A$ is weakly majorized by $B$ and one writes $A\prec_{w}B$. If furthemore in (a) one has the equality $\| A\|_{(n)}=\| B\|_{(n)}$, that it is $A$ and B have the same trace, then  $A$ is majorized by $B$, written $A\prec B$. Thus $A\prec B$ means that (c) holds with the equality sign: $A$ is in the convex hull of the unitary orbit of $B$. Theorem \ref{th2} is a  special majorization.  

A  linear map $\Phi:\bM_n\to \bM_n$ is called doubly stochastic if $\Phi$  preserves positivity, identity, and trace. For all $A\in\bM_n^+$, we then have
$\Phi(A)\prec A$, see the last section of Ando's survey \cite{A}.

\section{The distance from $0$ to the numerical range}

We state our main result and infer several corollaries. The proof of the theorem is postponed to Section 3. 

\vskip 10pt
\begin{theorem}\label{th-dist} Let $\begin{bmatrix} A &X \\
X^* &B\end{bmatrix} $ be a positive matrix  partitioned into four blocks
in $\bM_n$ and let $d={\mathrm{dist}}(0,W(X))$. Then, for all symmetric norms,
$$
 \left\| \begin{bmatrix} A &X \\
X^* &B\end{bmatrix}\right\| \ge
\left\| \left(\frac{A+B}{2} +dI\right) \oplus  \left(\frac{A+B}{2} -dI\right) \right\|.
$$
\end{theorem}

\vskip 5pt
Here, the direct sum  is a standard notation for   block-diagonal matrices
$$
X\oplus Y=\begin{bmatrix}X &0 \\ 0 &Y \end{bmatrix}.
$$

Since we have equality for the trace, Theorem \ref{th-dist} is a majorization relation. We  have $(A+B)/2 \ge dI$, otherwise, the trace norm of the left-hand side would be strictly smaller than the right-hand side one, a contradiction with the theorem.

 By a basic principle of majorization, Theorem \ref{th-dist} is equivalent to some trace inequalities.

\vskip 5pt
\begin{cor}\label{cor-tr}  Let $\begin{bmatrix} A &X \\
X^* &B\end{bmatrix} $ be a positive matrix  partitioned into four blocks
in $\bM_n$ and let $d={\mathrm{dist}}(0,W(X))$. Then, for every convex function $g:[0,\infty)\to(-\infty,\infty)$,
$$
{\mathrm{Tr\,}} g\left(\frac{A+B}{2} +dI\right) + {\mathrm{Tr\,}} g\left(\frac{A+B}{2} -dI\right)  \le
{\mathrm{Tr\,}} g\left(\begin{bmatrix} A &X \\
X^* &B\end{bmatrix}\right).
$$
\end{cor}

\vskip 5pt
Symmetric norms $\|\cdot\|$ on $\bM_n^+$ are the  homogeneous, unitarily invariant, convex functionals.  The concave counterpart, the symmetric anti-norms $\|\cdot\|_! $, have been introduced and studied in  papers \cite{BH1} and \cite[Section 4]{BH2}. We recall their basic properties, parallel to those of symmetric norms given at the end of Section 1. Symmetric anti-norms  on $\bM_n^+$
are continuous functionals characterized by  three properties:
\begin{itemize}
\item[(i)] $\| \lambda A\|_! = \lambda\|A\|_!$  for all $A\in \bM_n^+ $ and all  $\lambda\ge 0$,
\item[(ii)] $\| UAU^*\|_!$ for all $A\in \bM_n^+$ and all unitaries $U\in\bM_n$, 
\item[(iii)] $ \|A+B\|_! \ge   \|  A\|_!  + \|  B\|_!$ for all $A,\,B\in \bM_n^+ $.
\end{itemize}
 Let $\lambda_1^{\uparrow}(A)\le\cdots\le \lambda_n^{\uparrow}(A)$ stand for the eigenvalues of $A\in\bM_n^+$ arranged in non-decreasing order. Then, the Ky Fan $k$-anti-norms,
$$
\| A\|_{(k)!} =\sum_{j=1}^k \lambda_j^{\uparrow}(A)
$$
are symmetric anti-norms, $k=1,\ldots,n$.  The following conditions are equivalent:
\begin{itemize}
\item[(a)] $\| A\|_{(k)!} \ge \| B\|_{(k)!} $ for all $k=1,\ldots,n$,
\item[(b)] $\| A\|_! \ge \| B\|_!$ for all symmetric anti-norms,
\item[(c)] The vector of the eigenvalues of $A$ is dominated by some convex combination of permutations of the vector of the eigenvalues of $B$, equivalently, 
$$
A \ge \sum_{i=1}^{n+1} \alpha_i U_i BU_i^*
$$
for some  unitary matrices $U_i$ and some  scalars $\alpha_i\ge 0$ such that $\sum_{i=1}^{n+1}\alpha_i=1$.
\end{itemize}
 The
continuity assumption is not essential, but deleting it would lead to rather strange
 functionals which are not continuous on the boundary of $\bM_n^+$, such as
$\| A\|_!:={\mathrm{Tr\,}}A$ if $A$ is invertible and $\| A\|_!:=0$ if $A$ is not
invertible.

Note that the trace norm is both a symmetric norm and a symmetric anti-norm and that the majorization $A\prec B$ in $\bM_n^+$ also entails that $\| A\|_!\ge \| B\|_!$ for all symmetric anti-norms. Thus Theorem \ref{th-dist} is equivalent to the following statement:

\vskip 5pt
\begin{cor}\label{cor-anti} Let $\begin{bmatrix} A &X \\
X^* &B\end{bmatrix} $ be a positive matrix  partitioned into four blocks
in $\bM_n$, let $d={\mathrm{dist}}(0,W(X))$. Then, for all symmetric anti-norms,
$$
\left\| \left(\frac{A+B}{2} +dI\right) \oplus  \left(\frac{A+B}{2} -dI\right) \right\|_! \ge
\left\| \begin{bmatrix} A &X \\
X^* &B\end{bmatrix}\right\|_!.
$$
\end{cor}

\vskip 10pt
\begin{cor}\label{cor-maxmin} Let $\begin{bmatrix} A &X \\
X^* &B\end{bmatrix} $ be a positive matrix  partitioned into four blocks
in $\bM_n$ and let $d={\mathrm{dist}}(0,W(X))$. Then, 
$$
 \lambda_1^{\downarrow}\left(\begin{bmatrix} A &X \\
X^* &B\end{bmatrix}\right)-
 \lambda_1^{\downarrow}\left(\frac{A+B}{2}\right)  \ge d 
$$
and
$$
 \lambda_1^{\uparrow}\left(\frac{A+B}{2}\right) -
 \lambda_1^{\uparrow}\left(\begin{bmatrix} A &X \\
X^* &B\end{bmatrix}\right)\ge d.
$$
\end{cor}

\vskip 10pt
\begin{proof} The first inequality follows from Therorem \ref{th-dist} applied to the symmetric norm $A\mapsto\lambda_1^{\downarrow}(A)$ (the operator norm on the positive cone), while
 the second inequality follows from Corollary \ref{cor-anti} applied to the anti-norm 
 $A\mapsto\lambda_1^{\uparrow}(A)$
\end{proof}

\vskip 5pt
By adding these two inequalities we get an estimate for the spread of the matrices, i.e., for the diameter of the numerical ranges.

\vskip 10pt
\begin{cor}\label{cor-diam} For every  positive matrix  partitioned into four blocks of same size,
$$ 
{\mathrm{diam}}\, W\left(\begin{bmatrix} A &X \\
X^* &B\end{bmatrix}\right)
- {\mathrm{diam}}\, W\left(\frac{A+B}{2}\right) \ge 2d,
$$
where $d$ is the distance from $0$ to $W(X)$.
\end{cor}

\vskip 10pt
Of course
$$
{\mathrm{diam}}\, W\left(\begin{bmatrix} A &X \\
X^* &B\end{bmatrix}\right) \ge {\mathrm{diam}}\, W\left(\begin{bmatrix} A &0 \\
0 &B\end{bmatrix}\right)\ge  {\mathrm{diam}}\, W\left(\frac{A+B}{2}\right),
$$
however the ratio
$$
\rho=\frac{1}{2d}\left\{{\mathrm{diam}}\, W\left(\begin{bmatrix} A &X \\
X^* &B\end{bmatrix}\right) - {\mathrm{diam}}\, W\left(\begin{bmatrix} A &0 \\
0 &B\end{bmatrix}\right)\right\}
$$
can be arbitrarily small as shown by the following example where the blocks are in $\bM_2$,
$$
\begin{bmatrix} A &X \\
X^* &B\end{bmatrix}
=\begin{bmatrix} \begin{pmatrix} \alpha&0 \\ 0&\alpha^{-1}\end{pmatrix}& \begin{pmatrix} 1&0 \\ 0&1\end{pmatrix}\\
\begin{pmatrix} 1&0 \\ 0&1\end{pmatrix}&
\begin{pmatrix} \alpha^{-1}&0 \\ 0&\alpha\end{pmatrix}\end{bmatrix},
$$
and by noting that $\rho$ then takes the value $2/\alpha$ which tends to $0$ as 
$\alpha\to \infty$.

The Minkowki inequality for positive $m$-by-$m$ matrices,
$$
{\mathrm{det}}^{1/m}(A+B) \ge {\mathrm{det}}^{1/m}(A) + {\mathrm{det}}^{1/m}(B),
$$
shows that the functional $A\mapsto {\mathrm{det}}^{1/m}(A)$ is a symmetric anti-norm on $\bM_m^+$. For this anti-norm
Theorem \ref{th-dist} reads as:

\vskip 10pt
\begin{cor}\label{cor-det} Let $\begin{bmatrix} A &X \\
X^* &B\end{bmatrix} $ be a positive matrix  partitioned into four blocks
in $\bM_n$ and let $d={\mathrm{dist}}(0,W(X))$. Then, 
$$
\det \left\{\left(\frac{A+B}{2}\right)^2 -d^2I\right\} \ge
 \det \left(\begin{bmatrix} A &X \\
X^* &B\end{bmatrix}\right).$$
\end{cor}

\vskip 5pt
Letting $X=0$, we recapture a basic property: the determinant is a  log-concave map on the positive  cone of $\bM_n$. Hence Corollary \ref{cor-det} refines this property. 

By a basic principle of majorization, Corollary \ref{cor-anti} is equivalent to the following seemingly more general statement.

\vskip 10pt
\begin{cor}\label{cor-antif} Let $\begin{bmatrix} A &X \\
X^* &B\end{bmatrix} $ be a positive matrix  partitioned into four blocks
in $\bM_n$, let $d={\mathrm{dist}}(0,W(X))$, and let $f(t)$ be a nonnegative concave function on $[0,\infty)$. Then, 
$$
\left\| f\left(\frac{A+B}{2} +dI\right) \oplus  f\left(\frac{A+B}{2} -dI\right) \right\|_! \ge
\left\| f\left(\begin{bmatrix} A &X \\
X^* &B\end{bmatrix}\right)\right\|_!
$$
for all symmetric anti-norms.
\end{cor}

\section{Proof of Theorem \ref{th-dist}}

We want to show the majorization in $\bM_{2n}^+$
\begin{equation}\label{eqmaj}
\begin{bmatrix} \frac{A+B}{2}+dI &0 \\
0 &\frac{A+B}{2}-dI\end{bmatrix}
\prec
\begin{bmatrix} A &X \\
X^* &B\end{bmatrix}
\end{equation}
where $d={\mathrm{dist}(0,W(X)}$. We use two lemmas, the first one might belong to folklore.

\vskip 5pt
\begin{lemma}\label{lem1} Let $\{A_k\}_{k=1}^m$ and  $\{B_k\}_{k=1}^m$ be two families  of $r$-by-$r$ positive matrices such that $A_k\prec B_k$ for each $k$. Then,
$$
\oplus_{k=1}^m A_k \prec \oplus_{k=1}^m B_k.
$$
\end{lemma}

\vskip 5pt
\begin{proof} Let $p_k$ denote any integer such that $0\le p_k \le m$, $k=1,\ldots,m$.
With this notation, we then have, for each integer $p=1,\ldots,mr$, 
\begin{align*}
\sum_{j=1}^p\lambda_j^{\downarrow}\left(\oplus_{k=1}^m A_k \right)
&=\max_{p_1+p_2+\cdots+p_m=p}\,\sum_{k=1}^m\sum_{j=1}^{p_k}(A_k)\\
&\le \max_{p_1+p_2+\cdots+p_m=p}\,\sum_{k=1}^m\sum_{j=1}^{p_k}(B_k)\\
&=\sum_{j=1}^p\lambda_j^{\downarrow}\left(\oplus_{k=1}^m B_k \right)
\end{align*}
with equality for $p=mr$.
\end{proof}

\vskip 5pt
\begin{lemma}\label{lem2} Let $X,Y \in\bM_n^+$ and let $\delta>0$ be such that $X\ge Y\ge \delta I$. Then,
$$
\begin{bmatrix} X+\delta I&0 \\
0 &X-\delta I\end{bmatrix}
\prec
\begin{bmatrix} X+Y &0 \\
0&X-Y\end{bmatrix}.
$$
\end{lemma}

\vskip 5pt
\begin{proof} Let $\{e_k\}_{k=1}^n$ be an orthonormal basis of $\bC^n$ and define two $n$-by-$n$ diagonal positive matrices
$$
D_+=\mathrm{diag}(\langle e_1, (X+Y)e_1\rangle, \ldots, \langle e_n, (X+Y)e_n\rangle)
$$
and
$$
D_-=\mathrm{diag}(\langle e_1, (X-Y)e_1\rangle, \ldots, \langle e_n, (X-Y)e_n\rangle).
$$
Since extracting a diagonal is a doubly stochastic map (a pinching), we have
\begin{equation}\label{eq-DD}
\begin{bmatrix} D_+ &0 \\
0&D_-\end{bmatrix}\prec
\begin{bmatrix} X+Y &0 \\
0&X-Y\end{bmatrix}.
\end{equation}
Now, choose the basis  $\{e_k\}_{k=1}^n$ as a basis of eigenvectors for $X$, $\lambda_k^{\downarrow}(X)=\langle e_k, Xe_k\rangle$, and observe that the majorization in $\bM_2^+$,
$$
\begin{pmatrix}
  \lambda_k^{\downarrow}(X) + \delta& 0 \\ 0&  \lambda_k^{\downarrow}(X)-\delta
\end{pmatrix}
\prec
\begin{pmatrix}
 \langle e_k, (X+Y)e_k\rangle & 0 \\ 0&  \langle e_k, (X-Y)e_k\rangle 
\end{pmatrix},
$$
holds for every $k$. Applying Lemma \ref{lem1} then shows that
$$
\bigoplus_{k=1}^n \begin{pmatrix}
  \lambda_k^{\downarrow}(X) + \delta& 0 \\ 0&  \lambda_k^{\downarrow}(X)-\delta
\end{pmatrix}
\prec
\bigoplus_{k=1}^n
\begin{pmatrix}
 \langle e_k, (X+Y)e_k\rangle & 0 \\ 0&  \langle e_k, (X-Y)e_k\rangle
\end{pmatrix}.
$$
This means that
$$
\begin{bmatrix} X+\delta I&0 \\
0 &X-\delta I\end{bmatrix}
\prec
\begin{bmatrix}D_+ &0 \\
0&D_-\end{bmatrix}
$$
and we may combine this majorization with \eqref{eq-DD}  to complete the proof.
\end{proof}

\vskip 5pt
We turn to the proof of \eqref{eqmaj}.

\vskip 5pt
\begin{proof} 
 Suppose first that $d=0$, that is $0\in W(X)$. Note that
\begin{equation}\label{eq-1}
\begin{bmatrix} A &0 \\
0 &B\end{bmatrix}
\prec
\begin{bmatrix} A &X \\
X^* &B\end{bmatrix}
\end{equation}
as the operation of taking the block diagonal is doubly stochastic. 

Using the unitary congruence with
\begin{equation}\label{eqJ}
J=\frac{1}{\sqrt{2}}\begin{bmatrix} I &-I \\
I &I\end{bmatrix}
\end{equation}
 we observe that
\begin{equation*}
J\begin{bmatrix} A &0 \\
0 &B\end{bmatrix} J^*= \begin{bmatrix} \frac{A+B}{2} & \frac{A-B}{2}  \\
 \frac{A-B}{2} & \frac{A+B}{2} \end{bmatrix}
\end{equation*}
Hence we have 
\begin{equation*}
\begin{bmatrix} \frac{A+B}{2} & 0  \\
0 & \frac{A+B}{2} \end{bmatrix}\prec\begin{bmatrix} A &0 \\
0 &B\end{bmatrix}
\end{equation*}
and combining with $\eqref{eq-1}$ establishes \eqref{eqmaj} for the case $d=0$.

Now assume that $d>0$, that is $0\notin W(X)$. Using the unitary congruence implemented by
$$
\begin{bmatrix} I &0 \\
0 &e^{-i\theta}I\end{bmatrix}
$$
we may replace the right hand side of \eqref{eqmaj} with
$$
\begin{bmatrix} A &e^{i\theta}X \\
e^{-i\theta}X^* &B\end{bmatrix}
$$
Thanks to the rotation property $W(e^{i\theta}X)=e^{i\theta}W(X)$, by choosing the adequate $\theta$, we may then and do assume that $W(X)$ lies the half-plane of $\bC$ consiting of complex numbers with real parts greater or equal than $d$,
$$
W(X) \subset \{ z= x+iy \ : \ x\ge d\,\}.
$$
The projection property for the real part of the numerical range, ${\mathrm{Re}}\,W(X)=W({\mathrm{Re}}\,X)$ with ${\mathrm{Re}}\,X=(X+X^*)/2$, then ensures that
$$
{\mathrm{Re}}\,X \ge dI.
$$
Now, using again a unitary congruence with \eqref{eqJ}, wet get
\begin{equation*}
J\begin{bmatrix} A &X \\
X^* &B\end{bmatrix} J^*= \begin{bmatrix} \frac{A+B}{2} +{\mathrm{Re}} X& \ast \\
\ast & \frac{A+B}{2} -{\mathrm{Re}} X\end{bmatrix}
\end{equation*}
where $\ast$ stands for unspecified entries. Hence
\begin{equation*}
 \begin{bmatrix} \frac{A+B}{2} +{\mathrm{Re}} X& 0 \\
0 & \frac{A+B}{2} -{\mathrm{Re}} X\end{bmatrix} 
\prec
\begin{bmatrix} A & X \\
X^* &B\end{bmatrix}
\end{equation*}
and applying Lemma \ref{lem2} then yields \eqref{eqmaj}.
\end{proof}

\vskip 5pt

Jean-Christophe Bourin

Laboratoire de math\'ematiques de Besan\c{c}on, UMR n$^o$ 6623, CNRS,

Universit\'e de Bourgogne Franche-Comt\'e

Email: jcbourin@univ-fcomte.fr

  \vskip 15pt
Eun-Young Lee

 Department of mathematics, KNU-Center for Nonlinear
Dynamics,

Kyungpook National University,

 Daegu 702-701, Korea.

 Email: eylee89@knu.ac.kr


\begin{thebibliography}{99}
{\small

 \bibitem{A} T.\ Ando, Majorization, doubly stochastic matrices, and comparison of eigenvalues,  {\it Linear Algebra Appl.}\  118  (1989), 163--248. 






 \bibitem{Bhatia} R.\ Bhatia, Matrix Analysis, Gradutate Texts in Mathematics, Springer, New-York, 1996.

\bibitem{BH1} J.-C.\ Bourin and F.\ Hiai,
Norm and anti-norm inequalities for positive semi-definite matrices,
\textit{Internat.\ J.\ Math.}\ 22 (2011), 1121--1138.

\bibitem{BH2} J.-C.\ Bourin and F.\ Hiai, Jensen and Minkowski inequalities for
operator means and anti-norms, {\it Linear Algebra Appl.}\ 456 (2014), 22--53.


\bibitem{BL} J.-C. Bourin and E.-Y.\ Lee, Decomposition and partial trace of positive matrices with Hermitian blocks, Int. J. Math. 24 (2013) 1350010. 


\bibitem{BLL2} J.-C.\ Bourin, E.-Y.\ Lee and M.\ Lin, Positive matrices partitioned into a small number of Hermitian blocks,{\it Linear Algebra Appl}.\ 438 (2013) 2591--2598.

\bibitem{BM}
J.-C.\ Bourin and A.\ Mhanna, Positive block matrices and numerical ranges, {\it C.\ R.\ Acad.\ Sci.\ Paris},  355 no.\,10 (2017) 1077--1081. 


\bibitem{Cain} B.\ E.\ Cain, Improved inequalities for the numerical radius: when inverse commutes with the norm. Bull. Aust. Math. Soc.  97  (2018),  no. 2, 293--296.

\bibitem{GLRT} M.\ Gumus, J.\ Liu, S.\ Raouafi, T.-Y.\ Tam, Positive semi-definite $2\times2$  block matrices and norm inequalities, {\it Linear Algebra Appl}.\ 551 (2018) 83--91.

\bibitem{Hay} T.\ Hayashi, On a norm inequality for a positive block-matrix, Linear Algebra Appl.  566  (2019), 86--97. 

\bibitem{Hol} J.A.R.\ Holbrook, Multiplicative properties of the numerical radius in operator theory, J. Reine Angew. Math. 237 (1969), 166--174. 

\bibitem{HJ} R.\ Horn and C.R.\ Johnson Topics in matrix analysis. Corrected reprint of the 1991 original. Cambridge University Press, Cambridge, 1994.

\bibitem{K1} F.\ Kittaneh, Numerical radius inequalities for Hilbert space operators. Studia Math.  168  (2005),  no.\,1, 73--80. 

\bibitem{Mha} A.\ Mhanna, On symmetric norm inequalities and positive definite block-matrices, {\it Math.\ Ineq.\  Appl.}, 
  21  (2018),  no.\,1, 133–138.

\bibitem{S} B.\ Simon, 
Trace ideals and their applications. Second edition. Mathematical Surveys and Monographs, 120, American Mathematical Society, Providence, RI, 2005.

}

\end{thebibliography}
\end{document}